\documentclass[amsmath,12pt,a4paper]{amsart}
\usepackage{amsfonts, amssymb, amsmath, amscd}
\usepackage{amsthm}
\usepackage[pdftex]{graphicx}
\usepackage{multicol}

\usepackage{tikz}
\usepackage{ragged2e,caption}
\usepackage{bm,bigints, comment}
\usepackage{mathrsfs}
\usepackage{color}

\newtheorem{theorem}{Theorem}

\newtheorem{Corollary}{Corollary}

\newtheorem{lemma}{Lemma}

\begin{document}
\title[The Navier-Stokes Equation and Helmholtz decomposition]{The Navier-Stokes Equation and Helmholtz decomposition}
\author{Roy Burson}
\maketitle

\begin{abstract}
This work explores Navier-Stokes equation with no gravitational forces. In short it shows that any smooth solution that decays quickly must take the form 

$$ \textbf{u}(x,t)- \dfrac{1}{4\pi}\textbf{Curl}\Biggl( \int_{\mathbb{R}^3}^{}{\dfrac{\textbf{Curl} (\textbf{u} (x^\prime,t))}{|x-x^\prime|}}dV^\prime\Biggr) = -\int_{0}^{t}{\dfrac{1}{\rho}  \textbf{Grad}\big(\Gamma(x,s)\big)}ds.$$ 
 with $\Gamma$ a known function which is related to the heat equation. Consequently, any curl free solution must be written as
$$\textbf{u}(x,t) = -\dfrac{1}{\rho} \textbf{Grad}\biggl(\int_{0}^{t}{\Gamma(x,s) ds}\biggr).$$ Even further it shows if there exist a value $k\in \mathbb{N}$ such that $$\textbf{curl}^k\biggl((\textbf{u}\cdot \nabla )\textbf{u}\biggr)(x,t)=\textbf{0}$$ for all $0<t^\prime\le t$ then $$\textbf{u}(x,t) = \textbf{H}^{k+1}(\xi_1,\xi_2,\xi_3,t) -\int_{0}^{t}{\dfrac{1}{\rho}  \textbf{Grad}\big(\Gamma(x,s)\big)}ds, \indent t\in [t^\prime,\infty)$$ with $$\xi_i(x,t):= \int_{\mathbb{R}^3}^{}{\alpha(x-y,\dfrac{t}{\nu})v^k_i(x,0)}dy, \indent v^k_i(x,0) = \biggl(\textbf{curl}^k(\textbf{u}(x,0))\biggr)_i \indent 1\le i\le 3,$$ $\textbf{H}^k$ the $k^{th}$ application of Helmholtz operator, and $\alpha(x,t)$ the kernel to the heat equation. Thus if the representation found in this work is a solution and if the solution is unique, then this is the only possible solution. 
\end{abstract}

\section{The Navier-Stokes equation in $\mathbb{R}^3$}\label{S1}

The Navier-Stokes equation in $\mathbb{R}^3$ subjected to no gravitational forces are provided as:

\begin{align}
 &\dfrac{\partial }{\partial t} \textbf{u} + \left(\textbf{u}\cdot \nabla \right)\textbf{u} -    \nu \bigtriangleup \textbf{u} = -\dfrac{1}{\rho}  \textbf{Grad}(p)\label{Eq1}\\
&\textbf{Div}(\textbf{u})  = 0 \label{Eq2}\\
&\textbf{u}(x,0) = \textbf{u}^0(x).\label{Eq3}
\end{align}

Here $\rho>0$ is a scalar which represents the fluids density and $\nu>0$ is a scalar which denotes the fluids kinetic viscosity. These equations are to be solved for the velocity field $\textbf{u}(x,t)$ and scalar pressure function $p(x,t)$ as time progresses. For physical reasonable solutions the initial condition is restricted to be smooth (infinitely differentiable) and decay $$|\textbf{u}^0|_{L^2_{x}(\mathbb{R})}\le \infty$$ so that $\textbf{u}(x,t)$ does not grow large as $|x|\rightarrow \infty$. According to Fefferman's official problem description \cite{Feff2006} we say a solution $\textbf{u}$ and $p$ is reasonable if

\begin{equation}\label{Eq4}
p, \textbf{u}\in C^{\infty}(\mathbb{R}^3\times [0,\infty)) \quad \text{ and } \quad \int_{\mathbb{R}^3}^{}{|\textbf{u}|} < \epsilon, \quad\forall t >0.
\end{equation}

That is $\textbf{u}$ and $p$ are infinitely differentiable and $\textbf{u}$ has finite energy. 

\begin{theorem}\label{Thm1}
If there is a smooth solution to Eq.\eqref{Eq1}-Eq.\eqref{Eq3} that satisfy the conditions of Eq.\eqref{Eq4} such that $\left(\textbf{u}\cdot \nabla \right)\textbf{u} -    \nu \bigtriangleup \textbf{u}$ vanishes as fast as $\dfrac{1}{x}$ is twice differentiable, and is of bounded support then the velocity must take the form

$$ \textbf{u}(x,t)- \dfrac{1}{4\pi}\textbf{Curl}\Biggl( \int_{\mathbb{R}^3}^{}{\dfrac{\textbf{Curl} (\textbf{u} (x^\prime,t))}{|x-x^\prime|}}dV^\prime\Biggr) = -\int_{0}^{t}{\dfrac{1}{\rho}  \textbf{Grad}\big(\Gamma(x,s)\big)}ds, \indent t>0$$ with $$\Gamma(x,t):=\int_{\mathbb{R}^3}^{}{\alpha(x-y,\dfrac{t}{\nu})(p - \rho \phi)\circ(x,0)}dy$$ and $\alpha(x,t)$ the kernel to the non-homogenous heat equation with  $\phi(x,0)$ denoting the irrotational part of the non-linear term $\left(\textbf{u}\cdot \nabla \right)\textbf{u} -    \nu \bigtriangleup \textbf{u}$ at time $t=0$. 

\begin{proof}
 Assume for simplicity that there is a smooth solution to Eq.\eqref{Eq1}-Eq.\eqref{Eq3} that satisfy the conditions of Eq.\eqref{Eq4}. If the vector field $\left(\textbf{u}\cdot \nabla \right)\textbf{u} -    \nu \bigtriangleup \textbf{u}$ vanishes as fast as $\dfrac{1}{x}$ as $x\rightarrow \infty$, is twice differentiable, and is of bounded support then according to Helmholtz decomposition the vector field $\left(\textbf{u}\cdot \nabla \right)\textbf{u} -    \nu \bigtriangleup \textbf{u}$ can be written as 

\begin{equation}\label{Eq5}
\left(\textbf{u}\cdot \nabla \right)\textbf{u} -    \nu \bigtriangleup \textbf{u}= \textbf{Grad}(\phi) + \textbf{Curl}(\Phi)
\end{equation}

where $\phi$ and $\Phi$ are provided explicitly by the formula
 
$$
\phi(x) = \dfrac{1}{4\pi}\int_{\mathbb{R}^3}^{}{\dfrac{\textbf{Div}\bigl(\left(\textbf{u}\cdot \nabla \right)\textbf{u} -    \nu \bigtriangleup \textbf{u}\bigr)(x^\prime)}{|x-x^\prime|}}dV^\prime
$$ 

and 

$$ 
\Phi(x) = \dfrac{1}{4\pi}\int_{\mathbb{R}^3}^{}{\dfrac{\textbf{Curl}\bigl(\left(\textbf{u}\cdot \nabla \right)\textbf{u} -    \nu \bigtriangleup \textbf{u}\bigr)(x^\prime)}{|x-x^\prime|}}dV^\prime
$$ 

respectively (see \cite{Griffiths1999}). Substituting Eq.\eqref{Eq5} into Eq.\eqref{Eq1} yields a new expression for the velocity

\begin{equation}\label{Eq6}
\dfrac{\partial }{\partial t} \textbf{u}+\textbf{Curl}(\Phi) = -\dfrac{1}{\rho}  \textbf{Grad}(p - \rho \phi).
\end{equation}

Computing the divergence of Eq.\eqref{Eq6} one finds

\begin{equation}\label{Eq7}
\dfrac{\partial }{\partial t} \textbf{Div}(\textbf{u})+\textbf{Div}\big(\textbf{Curl}(\Phi)\big) = -\dfrac{1}{\rho}  \textbf{Div}\big(\textbf{Grad}(p - \rho \phi)\big).
\end{equation}

According to Eq.\eqref{Eq2} and the fact that $\textbf{Div}\big(\textbf{Curl}(\textbf{f})\big) = 0$ for any vector field $\textbf{f}$ Eq.\eqref{Eq7} then reduces to

\begin{equation}\label{Eq8}
\dfrac{1}{\rho}  \textbf{Div}\big(\textbf{Grad}(p - \rho \phi)\big)= \dfrac{1}{\rho}\sum_{i=1}^{3}{\dfrac{\partial^2}{\partial x_i^2} (p - \rho \phi)} = 0  .
\end{equation}

Conveniently, Eq.\eqref{Eq8} is simply the non-homogeneous heat equation in terms of the variable $(p - \rho \phi)$ which then yields the analytic result

\begin{equation}\label{Eq9}
(p - \rho \phi) = \int_{\mathbb{R}^3}^{}{\alpha(x-y,\dfrac{t}{\nu})(p - \rho \phi)\circ(x,0)}dy:= \Gamma(x,t)
\end{equation}

with $\alpha(x,t)$ the kernel to the non-homogeneous heat equation 

$$
\alpha(x,t):= \dfrac{1}{(4\pi t)^{3/2}}~ \exp\biggl(\dfrac{|x|^2}{4t}\biggr), \quad x\in \mathbb{R}^3,~~ t>0.$$

according to Evan's book on PDE's \cite[p.~48]{Evans2022} as the heat equations solution is unique. Next using Eq.\eqref{Eq9} we may rewrite Eq.\eqref{Eq6} in terms of the mapping $\Gamma$ as followed

\begin{equation}\label{Eq10}
\dfrac{\partial }{\partial t} \textbf{u}+\textbf{Curl}(\Phi) = -\dfrac{1}{\rho}  \textbf{Grad}\big(\Gamma(x,t)\big).
\end{equation}

Recall $\Gamma(x,t)$ is known and the vector field $\Phi$ is not because it is a function of the velocity $\textbf{u}(x,t)$. Substituting the definition of $\Phi$ into Eq.\eqref{Eq10} reveals
 
\begin{equation}\label{Eq11}
\begin{split}
\dfrac{\partial }{\partial t} \textbf{u} + \dfrac{1}{4\pi}\textbf{Curl}&\left(\int_{\mathbb{R}^3}^{}{\dfrac{ \textbf{Curl}\bigl((\textbf{u}\cdot \nabla )\textbf{u} - \nu \bigtriangleup \textbf{u}\bigr)(x^\prime,t)}{|x-x^\prime|}}dV^\prime \right) \\
&\qquad= -\dfrac{1}{\rho}  \textbf{Grad}\big(\Gamma(x,t)).
\end{split}
\end{equation}

The exceptional part in this formula is that it does not involve any pressure even as difficult as it may seem.  A simple reduction can be done by comparing to the curl of Eq.\eqref{Eq1} and noticing the curl of the vector field $\textbf{u}(x,t)$ is equivalent to the curl of its non-linear part minus the Laplacian term

\begin{equation}\label{Eq12}
- \dfrac{d}{dt}\textbf{Curl}(\textbf{u}) = \textbf{Curl}\bigl(\left(\textbf{u}\cdot \nabla \right)\textbf{u} -\nu \bigtriangleup \textbf{u}\bigr).
\end{equation}

Substituting Eq.\eqref{Eq12} into Eq.\eqref{Eq11} gives the result

\begin{equation}\label{Eq13}
\dfrac{\partial }{\partial t} \textbf{u}(x,t)- \dfrac{1}{4\pi}\textbf{Curl}\Biggl(\int_{\mathbb{R}^3}^{}{ \dfrac{ \textbf{Curl}(\frac{d}{dt}\textbf{u} (x^\prime,t))}{|x-x^\prime|}}dV^\prime\Biggr) = -\dfrac{1}{\rho}  \textbf{Grad}\big(\Gamma(x,t)\big).
\end{equation}

We can pull out the time derivative (as the volume integral is independent of the time variable $t$) and we should actually get the following

\begin{equation}\label{Eq14}
 \textbf{u}(x,t)- \dfrac{1}{4\pi}\textbf{Curl}\Biggl( \int_{\mathbb{R}^3}^{}{\dfrac{\textbf{Curl} (\textbf{u} (x^\prime,t))}{|x-x^\prime|}}dV^\prime\Biggr) = -\int_{0}^{t}{\dfrac{1}{\rho}  \textbf{Grad}\big(\Gamma(x,s)\big)}ds.
\end{equation}

which proves the result.

\end{proof}
\end{theorem}

\begin{lemma}\label{Lem1}
If $\textbf{curl}(\textbf{u}(x,t^\prime))= \textbf{0}$ then  $\textbf{curl}(\textbf{u}(x,t)) = \textbf{0}$ for all $t\ge t^\prime>0$.

\begin{proof}
The proof is simple and instructive. Note that if $\textbf{u}(x,t)$ is smooth and satisfied Eq.\eqref{Eq1} then $\textbf{v}(x,t):= \textbf{curl}(\textbf{u}(x,t))$ satisfies $$
\dfrac{\partial}{\partial t}\textbf{v} + \textbf{curl}\biggl(\left(\textbf{u}\cdot \nabla \right)\textbf{u}\biggr) = \nu \bigtriangleup\textbf{v}.
$$ Next assume $\textbf{v}(x,t^\prime)= \textbf{0}$ for some $t^\prime>0$ then it follows that

\begin{equation}\nonumber
\begin{split}
\textbf{curl}\biggl(\left(\textbf{u}\cdot \nabla \right)\textbf{u}\biggr)(x,t^\prime) & = \biggl((\textbf{v}\cdot \nabla )\textbf{u} - (\textbf{u}\cdot \nabla )\textbf{v} \biggr)(x,t^\prime)=0 \textbf{0},
\end{split}
\end{equation}

and $\nu \bigtriangleup \textbf{v} (x,t^\prime) = \textbf{0}$ which implies $\dfrac{\partial}{\partial t}\textbf{v} (x,t^\prime) = \textbf{0}$. However, both $\dfrac{\partial}{\partial t}\textbf{v} (x,t^\prime) = \textbf{0}$ and $\textbf{v}(x,t^\prime)= \textbf{0}$ imply $\textbf{v}(x,t) = \textbf{0}$ for all $t\ge t^\prime>0$. 
\end{proof}
\end{lemma}

\begin{Corollary}\label{Cor1}
If there is a smooth solution to Eq.\eqref{Eq1}-Eq.\eqref{Eq3} that satisfy the conditions of Eq.\eqref{Eq4} such that $\left(\textbf{u}\cdot \nabla \right)\textbf{u} -    \nu \bigtriangleup \textbf{u}$ vanishes as fast as $\dfrac{1}{x}$ is twice differentiable and is of bounded support with $\textbf{curl}(\textbf{u}(x,t^\prime))= \textbf{0}$ for some $t^\prime \in (0,\infty)$ then $\textbf{u}$ must take the form  $$
\textbf{u}(x,t) = -\dfrac{1}{\rho} \textbf{Grad}\biggl(\int_{0}^{t}{\Gamma(x,s) ds}\biggr), \indent t\ge t^\prime>0.
$$ where $\Gamma$ defined explicitly in Eq.\eqref{Eq9}.

\begin{proof}
Assume the conditions of the hypothesis. By Lemma \ref{Lem1} if $\textbf{curl}(\textbf{u}(x,t^\prime))= \textbf{0}$ then  $\textbf{curl}(\textbf{u}(x,t))$ for all $t\ge t^\prime$. As an immediate consequence of Theorem \ref{Thm1} one must have 

$$
\textbf{u}(x,t) = -\dfrac{1}{\rho} \textbf{Grad}\biggl(\int_{0}^{t}{\Gamma(x,s) ds}\biggr), \indent t\ge t^\prime>0.
$$ 

\end{proof} 
\end{Corollary}

\begin{theorem}\label{Thm2}
If there is a smooth solution to Eq.\eqref{Eq1}- Eq.\eqref{Eq4} such that the non-linear component $\left(\textbf{u}\cdot \nabla \right)\textbf{u} -    \nu \bigtriangleup \textbf{u}$ vanishes as fast as $\dfrac{1}{x}$, is twice differentiable, has bounded support, and there exist a value $k\in \mathbb{N}$ such that $$\textbf{curl}^k\biggl((\textbf{u}\cdot \nabla)\textbf{u}\biggr)=\textbf{0}, \indent t\ge t^\prime$$ then the velocity $\textbf{u}(x,t)$ can be written as $$\textbf{u}(x,t) = \textbf{H}^{k+1}(\xi_1,\xi_2,\xi_3,t) -\int_{0}^{t}{\dfrac{1}{\rho}  \textbf{Grad}\big(\Gamma(x,s)\big)}ds$$ for $t\ge t^\prime>0$ with $$\xi_i(x,t):= \int_{\mathbb{R}^3}^{}{\alpha(x-y,\dfrac{t}{\nu})v^k_i(x,0)}dy, \indent v^k_i(x,0) = \biggl(\textbf{curl}^k(\textbf{u}(x,0))\biggr)_i, \indent 1\le i\le 3$$ and $\Gamma$ defined explicitly in Eq.\eqref{Eq9}, $\textbf{v}^k(x,t)$ the $k^{th}$ application of the curl operator

\begin{equation}\nonumber
\textbf{v}^k(x,t) = \textbf{curl}^k(\textbf{u}(x,t)) = \underbrace{\textbf{curl}(\textbf{curl}\cdots \cdots \textbf{curl})(\textbf{u}(x,t))}_{k-applications},
\end{equation} 

and $\textbf{H}^k(\textbf{v})$ defined as the kth composition of Helmholtz operator 

\begin{equation}\nonumber
\textbf{H}^k:\mathbb{R}^3\times [0,t)\rightarrow \mathbb{R}^3 \times [0,t), \indent \textbf{H}^k(\textbf{v}) := \underbrace{\biggl(\textbf{H}\circ \textbf{H}\circ \cdots \circ \textbf{H}\biggr) \circ (\textbf{v})}_{k-applications} \indent  \textbf{H}^1(\textbf{v})  = \dfrac{1}{4\pi} \int_{\mathbb{R}^3}^{}{\dfrac{\textbf{curl}(\textbf{v}(x^\prime,t))}{|x-x^\prime|}}dV^\prime.
\end{equation}

\begin{proof}
Assume the conditions of the hypothesis.  According to Eq.(\ref{Eq1}) the $k^{th}$ curl must satisfy the formula

\begin{equation}\label{Eq15}
\dfrac{\partial}{\partial t} \textbf{v}^k + \textbf{curl}^k\biggl(\left(\textbf{u}\cdot \nabla \right)\textbf{u}\biggr) = \nu \bigtriangleup \textbf{v}^k.
\end{equation}

Since $\textbf{curl}^k\biggl(\left(\textbf{u}\cdot \nabla \right)\textbf{u}\biggr) = \textbf{0}$  for all $t\ge t^\prime$ then immediately

\begin{equation}\label{Eq16}
\dfrac{\partial}{\partial t} \textbf{v}_i^k  = \nu \bigtriangleup \textbf{v}_i^k\indent 1\le i \le 3, \indent t\ge t^\prime>0.
\end{equation} 

and so $v_i^k$ admits a solution to the non-homogeneous heat equation

\begin{equation}\label{Eq17}
v_i^k(x,t)  = \int_{\mathbb{R}^3}^{}{\alpha(x-y,\dfrac{t}{\nu})v^k_i(x,0)}dy:= \xi_i(x,t), \indent v^k_i(x,0) = \biggl(\textbf{curl}^k(\textbf{u}(x,0))\biggr)_i, \indent 1\le i\le 3.
\end{equation}

By Helmholtz theorem as long as $\textbf{curl}^k(\textbf{u}(x,t))$ decays as fast as $\dfrac{1}{x}$ as $x\rightarrow \infty$  (which is guaranteed because $\textbf{v}^1$ is assumed to decay this quickly) then at the value $k$ we may write

\begin{equation}\label{Eq18}
\begin{split}
\textbf{v}(x,t)  & = \dfrac{1}{4\pi} \int_{\mathbb{R}^3}^{}{\dfrac{\textbf{curl}(\textbf{v}(x^\prime,t))}{|x-x^\prime|}}dV^\prime. \\\\
&=  \textbf{H}(\textbf{curl}(\textbf{v}(x,t)))\\
& = \underbrace{\biggl(\textbf{H}\circ \textbf{H}\circ \cdots \circ \textbf{H}\biggr)}_{l-times}\circ(\textbf{curl}^l(\textbf{v}(x,t),t), \indent L \le k.\\
& = \underbrace{\biggl(\textbf{H}\circ \textbf{H}\circ \cdots \circ \textbf{H}\biggr)}_{k-times}\circ(\xi_1(x,t),\xi_2(x,t),\xi_3(x,t),t).\\
\end{split}
\end{equation} 

This is simply a iteration of Helmholtz theorem on each $\textbf{v}^\epsilon$ until it terminates at $k$. Recall each $\xi_i$ is a function of space and time $(x,t)$ and thus this composition makes since. In this case using Theorem \ref{Thm1} we may write Eq.(\ref{Eq14}) as 

\begin{equation}\label{Eq19}
 \textbf{u}(x,t) = \textbf{H}^{k+1}(\xi_1,\xi_2,\xi_3,t) -\int_{0}^{t}{\dfrac{1}{\rho}  \textbf{Grad}\big(\Gamma(x,s)\big)}ds, \indent t\ge t^\prime.
\end{equation}

\end{proof}
\end{theorem}

\section{Relative Boundaries on the Vorticity Equation}
In the previous section we have shown that a solution to the Navier-stokes equation depends entirely on the vorticity whenever a smooth solution to the vorticity exist (that decays quickly) as we can write the velocity using Helmholtz decomposition. This is why we seek an analytical solution for the vorticity formula. As turbulence and vorticity is difficult to study (at least analytically) we seek another approach to examine the behavior. One probable approach is to bound the vorticity with some PDE we can solve. If we cannot bound the vorticity with a solvable PDE then it might actually contain a blow up. The intuition here is that we should be able to find a vector field that has similar speed to the vorticity but is simpler to solve with the same initial condition. We beg the following question: For each $1\le i\le 3$ is there a smooth mapping $F$ such that the PDE

\begin{equation}\label{Eq20}
\dfrac{\partial}{\partial t}\phi_i + F(\textbf{curl}(\Phi), \Phi, \nabla \phi_i, \nabla (\textbf{curl}(\Phi)_i)) = \nu\bigtriangleup \phi_i 
\end{equation}

has a solution $\Phi = [\phi_1,\phi_2,\phi_3]$ (which doesn't need to be unique) with

\begin{equation}\label{Eq21}
\dfrac{\partial}{\partial t}v_i \le \dfrac{\partial}{\partial t}\phi_i,\indent \phi_i(x,0) = v_i(x,0), \indent \phi_i(x,t)\le M  \indent 1\le i\le 3.
\end{equation}

That is the initial condition of $\phi_i$ matches $v_i$ and $\phi_i$ is bounded for each $i$. If no such $F$ exist then we might speculate that a blow up occurs for some initial condition. The first attempt one should take is to bound the no-linear term inside the vorticity equation.

\begin{theorem}\label{Thm3}
If there is a solution to Eq.(\ref{Eq1})- Eq.(\ref{Eq3})  such that $\nu \bigtriangleup v_i \ge 0 $ then  $$
\dfrac{1}{4}\dfrac{\partial}{\partial t} v_i \le |\textbf{v}|^2+ |\nabla v_i |^2 + |\textbf{u}|^2+  |\nabla u_i|^2 + \nu \bigtriangleup v_i.
$$ for all points in time $t>0$.

\begin{proof}
Define $\alpha_i$ as the angle between $\textbf{v}$ and $\nabla u_i$, and $\beta_i$ the angle between the vectors $\textbf{u}$ and $\nabla v_i$. Using the dot product identity the vorticity equation can be written as

\begin{equation}\label{Eq22}
\begin{split}
\dfrac{\partial}{\partial t} v_i &  = \biggl[|\textbf{v}|~|\nabla u_i|\cos(\beta_i)- |\textbf{u}|~|\nabla v_i|\cos(\alpha_i)\biggr]+  \nu \bigtriangleup v_i.\\
\end{split}
\end{equation}

Notice that if $\bigtriangleup v_i\ge 0$ for $1\le i\le 3$ then 

\begin{equation}\label{Eq23}
\begin{split}
\dfrac{\partial}{\partial t} v_i &  = \biggl[|\textbf{v}|~|\nabla u_i|\cos(\beta_i)- |\textbf{u}|~|\nabla v_i|\cos(\alpha_i)\biggr]+  \nu \bigtriangleup v_i\\
& \le  \biggl[|\textbf{v}|~|\nabla u_i|+ |\textbf{u}|~|\nabla v_i|\biggr]+  \nu \bigtriangleup v_i \\
& \le \max(|\textbf{v}|, |\nabla u_i |)^2 + \max(|\textbf{u}|, |\nabla v_i |)^2 +\nu \bigtriangleup v_i\\
& \le \biggl(|\textbf{v}|^2 + |\nabla u_i|^2\biggr) + \biggl(|\textbf{u}|^2 + |\nabla v_i|^2\biggr) + \nu \bigtriangleup v_i\\
&\le  2\max\biggl(|\textbf{v}|^2, |\textbf{u}|^2\biggr) +  2\max\biggl(|\nabla v_i|^2, |\nabla u_i|^2\biggr) +\nu \bigtriangleup v_i
\end{split}
\end{equation}

Thus, whenever $\bigtriangleup v_i\ge 0$ at least one of the four cases must hold at each point in time

\begin{equation}\label{Eq24}
\begin{split}
(i)~~~~\dfrac{\partial}{\partial t} v_i\le 2 \biggl(|\textbf{v}|^2 + |\nabla v_i|^2  \biggr)+ \nu \bigtriangleup v_i \\
\\
(ii)~~~~\dfrac{\partial}{\partial t} v_i\le 2 \biggl(|\textbf{u}|^2 + |\nabla u_i|^2  \biggr)+ \nu \bigtriangleup v_i \\
\\
(iii)~~~~\dfrac{\partial}{\partial t} v_i\le 2 \biggl(|\textbf{v}|^2 + |\nabla u_i|^2  \biggr)+ \nu \bigtriangleup v_i \\
\\
(iv)~~~~\dfrac{\partial}{\partial t} v_i\le 2 \biggl(|\textbf{u}|^2 + |\nabla v_i|^2  \biggr)+ \nu \bigtriangleup v_i. \\
\end{split}
\end{equation}

In any situation we may add (i)-(iv) to find

\begin{equation}\label{Eq25}
\dfrac{1}{4}\dfrac{\partial}{\partial t} v_i \le |\textbf{v}|^2+ |\nabla v_i |^2 + |\textbf{u}|^2+  |\nabla u_i|^2 + \nu \bigtriangleup v_i
\end{equation}

\end{proof}
\end{theorem}

\begin{theorem}\label{Thm4}
If there is a solution to Eq.(\ref{Eq1})- Eq.(\ref{Eq3})  such that $\nu \bigtriangleup v_i < 0 $ then  $$
\dfrac{1}{4}\dfrac{\partial}{\partial t} v_i \le |\textbf{v}|^2+ |\nabla v_i |^2 + |\textbf{u}|^2+  |\nabla u_i|^2.
$$ for all points in time $t>0$.
\begin{proof}
The proof is the same as Theorem \ref{Thm3}. Note that if $\nu \bigtriangleup v_i < 0 $ then 

\begin{equation}\label{Eq26}
\begin{split}
\dfrac{\partial}{\partial t} v_i &  = \biggl[|\textbf{v}|~|\nabla u_i|\cos(\beta_i)- |\textbf{u}|~|\nabla v_i|\cos(\alpha_i)\biggr]+  \nu \bigtriangleup v_i\\
& \le  \biggl[|\textbf{v}|~|\nabla u_i|+ |\textbf{u}|~|\nabla v_i|\biggr]\\
& \le \max(|\textbf{v}|, |\nabla u_i |)^2 + \max(|\textbf{u}|, |\nabla v_i |)^2 \\
& \le \biggl(|\textbf{v}|^2 + |\nabla u_i|^2\biggr) + \biggl(|\textbf{u}|^2 + |\nabla v_i|^2\biggr)\\
&\le  2\max\biggl(|\textbf{v}|^2, |\textbf{u}|^2\biggr) +  2\max\biggl(|\nabla v_i|^2, |\nabla u_i|^2\biggr) 
\end{split}
\end{equation}

Thus, whenever $\bigtriangleup v_i< 0$ at least one of the four cases must hold at each point in time

\begin{equation}\label{Eq27}
\begin{split}
(i)~~~~\dfrac{\partial}{\partial t} v_i\le 2 \biggl(|\textbf{v}|^2 + |\nabla v_i|^2  \biggr) \\
\\
(ii)~~~~\dfrac{\partial}{\partial t} v_i\le 2 \biggl(|\textbf{u}|^2 + |\nabla u_i|^2  \biggr)\\
\\
(iii)~~~~\dfrac{\partial}{\partial t} v_i\le 2 \biggl(|\textbf{v}|^2 + |\nabla u_i|^2  \biggr) \\
\\
(iv)~~~~\dfrac{\partial}{\partial t} v_i\le 2 \biggl(|\textbf{u}|^2 + |\nabla v_i|^2  \biggr). \\
\end{split}
\end{equation}

In any situation we may add (i)-(iv) to find 
 
 \begin{equation}\label{Eq28}
\dfrac{1}{4}\dfrac{\partial}{\partial t} v_i \le |\textbf{v}|^2+ |\nabla v_i |^2 + |\textbf{u}|^2+  |\nabla u_i|^2.
 \end{equation}
which completes the proof.\\
\end{proof}
\end{theorem}

\textbf{Note}: Theorem \ref{Thm3} and \ref{Thm4} illustrate that the time derivative of the vorticity is bounded by norm of four distinct vector fields squared ( in particular $|\textbf{v}|^2$, $|\nabla v_i |^2$, $|\textbf{u}|^2$, and $|\nabla u_i|^2$). As long as the Laplacian is non zero $\bigtriangleup v_i \neq 0$ there should exist non zero mappings $g_1$, $g_2$, $g_3$, and $g_4$ such that 

\begin{equation}\label{Eq29}
\begin{split}
 |\textbf{v}|^2\le g_1 \nu \bigtriangleup v_i, \indent  |\nabla v_i |^2\le g_2 \nu \bigtriangleup v_i,\indent  |\textbf{u}|^2\le g_3 \nu \bigtriangleup v_i, \indent  |\nabla u_i|^2\le g_4 \nu \bigtriangleup v_i
\end{split}
\end{equation}

since the four norms $|\textbf{v}|^2$, $|\nabla v_i |^2$, $|\textbf{u}|^2$, and $|\nabla u_i|^2$  are all real numbers in which case according to Eq.(\ref{Eq25}) we have

\begin{equation}\label{Eq30}
\dfrac{1}{4}\dfrac{\partial}{\partial t} v_i\le \biggl(g_1+g_2+g_3+g_4+1\biggr)\nu \bigtriangleup v_i.
\end{equation}

Thus hinting at a possible formulation for $F$ leading us to the following PDE

\begin{equation}\label{Eq31}
\dfrac{\partial}{\partial t} \lambda_i= \delta(x,t) \bigtriangleup \lambda_i
\end{equation}

with initial condition 

$$\lambda_i(x,0) = v_i(x,0) = (\textbf{curl}(\textbf{u}(x,0)))_i.$$

We wish to study the possibility of a blow up relative to $\delta$. In order to discover a possible boundary here we need to choose the correct delta so that all of this computation makes since. That is we must choose $$\delta = (g_1 +g_2 +g_3+g_4 +1)$$ for any $g_1,g_2,g_3$ and $g_4$ such that the inequalities in Eq.(\ref{Eq29}) hold valid for each point in time. For each $i$ one possible choice is $$g_1 = \dfrac{|\textbf{v}|^2}{\nu \bigtriangleup v_i},$$ $$ g_2 = \dfrac{|\nabla v_i|^2}{\nu \bigtriangleup v_i},$$   $$ g_3 = \dfrac{|\textbf{u}|^2}{\nu \bigtriangleup v_i},$$  $$ g_4 = \dfrac{|\nabla u_i|^2}{\nu \bigtriangleup v_i}.$$

If each $g_i$ does not blow in finite time then it is reasonable to speculate that $v_i$ will not blow up since the PDE in Eq.(\ref{Eq31}) should not blow up when $\delta$ decays and the solution to this particular PDE should always be larger then the vorticity because it starts at the same value but can have a larger derivative for all time. However, keeping track of all four of these quantities seems a bit challenging. Moreover, we can only use these functions whenever the Laplacian is non-zero. If the Laplacian is zero we must use a different method. We are interested in the limit 

\begin{equation}\label{Eq32}
\lim_{t\rightarrow \infty}{\delta(x,t)} = \lim_{t\rightarrow \infty}{\dfrac{|\textbf{v}|^2+ |\textbf{u}|^2+|\nabla v_i|^2+|\nabla u_i|^2}{\nu \bigtriangleup v_i}}+1.
\end{equation} \\

Whenever the Laplacian is positive and decays slower then the norms of $\textbf{v}$, $\textbf{u}$, $\nabla u_i$, and $\nabla v_i$ we must have

\begin{equation}\label{Eq33}
\lim_{t\rightarrow \infty}{\delta(x,t)} = 1
\end{equation} \\

This happens whenever $g_1$, $g_2$, $g_3$, and $g_4$ decay. We do not have a analytical formula for each $g_i$ but we know it should exist whenever a solution exist because they are just real numbers. In which case the vorticity seems to be bounded by a scaled version of the heat equation. Otherwise, a blow up is still possible. This is the reason that keeping track of the Laplacian of vorticity $\bigtriangleup v_i$ is important as it effects how $\delta$ can grow. When the Laplacian is zero then the derivative of the vorticity is bounded by the sum of norms of these four vector fields squared. So we need to examine each case separately or find out when the Laplacian is positive. If we can show that the quantities $|\textbf{v}|^2$, $|\textbf{u}|^2$, $|\nabla u_i|^2$, and $\nabla v_i$ decay faster then the denominator $\nu\bigtriangleup v_i$ then each $g_i$ will eventually approach zero without ever developing a singularity. It seems that whenever a solution decays the action of computing the norm should decrease the value more then taking the Laplacian. 

\section{Summary}
We take a second to reflect some of the results we found thus far. In the first section we have shown that under certain decay assumptions if the $k^{th}$ application of the curl operator of the non-linear term ever vanishes then the solution of the Selondial part of $\textbf{u}(x,t)$ can be written as a finite compositions of Helmholtz operator and a position vector related to the heat equation. In full generally, either $\textbf{curl}^k((\textbf{u}\cdot \nabla)\textbf{u}) = \textbf{0}$ for some $k$ or else $\textbf{curl}^k((\textbf{u}\cdot \nabla)\textbf{u})$ is infintly curlable. There is no other possibility.  This shows that the solution directly depends on its curl (whenever a smooth solution exist not necessarily a solution that blows up). We have yet to show there is a solution with $\textbf{curl}^k((\textbf{u}\cdot \nabla)\textbf{u}) \neq \textbf{0}$  for all $k\in \mathbb{N}$. If Eq.\eqref{Eq19} is a solution to Eq.\eqref{Eq1} and if there is another solution that is infintly curlable then the solution is not unique. Otherwise, if Eq.\eqref{Eq19} is a solution and the solution is unique then this is the only solution. Therefore, uniqueness of the solution would play a critical role in determining which representation the solution takes. In fact this work never showed that this formulation was a solution it only shows that a solution must take this form under these constraints (which in fact are not terrible to assume). It is still probable that a solution exist that has a blow up where one can compute the curl of the non-linear term consecutively.\\

\newpage
\bibliographystyle{amsplain}

\end{document}